\documentclass[arXiv]{actapoly}
\usepackage{subcaption}
\usepackage[]{algorithm2e}
\usepackage{bm}
\usepackage[thinc]{esdiff}
\usepackage{stanli}
\usepackage{multirow}
\usepackage{pgfplots}
\usepackage{dblfloatfix}
\pgfplotsset{width=8cm,compat=1.9}
\usetikzlibrary{arrows.meta,plotmarks}

\newcommand*\circled[1]{\tikz[baseline=(char.base)]{
		\node[shape=circle,draw,inner sep=0pt,fill=white, minimum size=4mm] (char) {#1};}}
\newcommand*\squared[1]{\tikz[baseline=(char.base)]{
		\node[shape=rectangle,draw,inner sep=0pt, minimum size=4mm] (char) {#1};}}

\begin{document}

\title{On optimum design of frame structures}

\correspondingauthor[M. Tyburec]{Marek Tyburec}{my}{marek.tyburec@fsv.cvut.cz}
\author[J. Zeman]{Jan Zeman}{my,their1}
\author[M. Kru\v{z}\'{i}k]{Martin Kru\v{z}\'{i}k}{their,their1}
\author[D. Henrion]{Didier Henrion}{their2,their3}

\institution{my}{Departement of Mechanics, Faculty of Civil Engineering, Czech Technical University in Prague, Th\'{a}kurova 7, 166 29 Prague 6, Czech Republic}
\institution{their}{Departement of Physics, Faculty of Civil Engineering, Czech Technical University in Prague, Th\'{a}kurova 7, 166 29 Prague 6, Czech Republic}
\institution{their1}{Institute of Information Theory and Automation, Czech Academy of Sciences, Pod vod\'{a}renskou v\v{e}\v{z}\'{i} 4, 182 08 Prague 8, Czech Republic }
\institution{their2}{LAAS-CNRS, University of Toulouse, 7 Avenue du Colonel Roche,  31400 Toulouse, France}
\institution{their3}{Department of Control Engineering, Faculty of Electrical Engineering, Czech Technical University in Prague, Karlovo n\'{a}m\v{e}st\'{i} 13, 121 35 Prague 2, Czech Republic}

\begin{abstract}
Optimization of frame structures is formulated as a~non-convex optimization problem, which is currently solved to local optimality. In this contribution, we investigate four optimization approaches: (i) general non-linear optimization, (ii) optimality criteria method, (iii) non-linear semidefinite programming, and (iv) polynomial optimization. We show that polynomial optimization solves the frame structure optimization to global optimality by building the (moment-sums-of-squares) hierarchy of convex linear semidefinite programming problems, and it also provides guaranteed lower and upper bounds on optimal design. Finally, we solve three sample optimization problems and conclude that the local optimization approaches may indeed converge to local optima, without any solution quality measure, or even to infeasible points. These issues are readily overcome by using polynomial optimization, which exhibits a finite convergence, at the prize of higher computational demands.
\end{abstract}


\maketitle

\section{Introduction}

Designing economical, efficient, and sustainable structures represents a major challenge of the contemporary society. While structural engineers literally explore designs that must satisfy the requirements of limit states analysis, there is usually an infinite number of such designs. Regardless of the properties of this \textit{feasible design space}, it is required to select only one, the \textit{best} design. This design quality is measured by an objective function, which usually approximates the expenses of production. Greatly simplified, one such (most common) criterion considers maximization of structural stiffness (while the amount of available material is limited), or, equivalently, minimization of structural volume (while requiring a certain structural stiffness) \cite{Bendsoe2003}.

Among the structural optimization problems the greatest progress has been achieved so-far in optimizing (the cross-sectional areas of) truss structures~\cite{Dorn1964}. This development can be attributed to convexity of the feasible design space~\cite[Sections 1.3.5, 3.4.3, and 4.8]{Ben-Tal2001}, as the (axial) structural stiffness is an affine function of the cross-sectional areas. On the other hand, when the \textit{bending} stiffness comes into a~consideration, convexity does not hold in general. Quite surprisingly, these problems are also much less studied, and to our knowledge only local optimization algorithms have yet been developed \cite{Steven2000,Rozvany1989}.

In this contribution, we investigate the problem of optimum design of frame structures. In particular, we develop four different methods in Section~\ref{sec:opt}: (i)~general non-linear optimization solved by the interior-point method of \texttt{fmincon}, (ii) the first-order optimality criteria (OC) method~\cite{Bendsoe2003,Rozvany1989}, (iii) reformulation of the problem into a non-linear semidefinite program (NSDP) solved by \textsc{PenLab}~\cite{Fiala2013}, and (iv) a suitably modified polynomial optmization (PO) problem (iii) solved globally using polynomial optimization methods~\cite{Lasserre2015} and the \textsc{Mosek}~\cite{mosek} optimizer. We show that the latter PO approach generates guaranteed lower and upper bounds on the solution, providing a means of assessing the design quality. Finally, Section~\ref{sec:sample} introduces a set of three sample optimization problems to compare the optimization approaches.

\section{Solution techniques to frame optimization}\label{sec:opt}

\subsection{Problem statement}

\addtocounter{equation}{2}
\begin{figure*}[t]
	\begin{equation}\label{eq:stiffmat}
	\renewcommand*{\arraystretch}{1.2}
	\mathbf{K}_i (a_i) = 
	\begin{pmatrix}
	\frac{E_i a_i}{\ell_i} & 0 & 0 & -\frac{E_i a_i}{\ell_i} & 0 & 0\\
	
	& \frac{12 E_i I_i}{\ell_i^3} & \frac{6 E_i I_i}{\ell_i^2} & 0 & -\frac{12 E_i I_i}{\ell_i^3} & \frac{6 E_i I_i}{\ell_i^2}\\
	
	&  & \frac{4 E_i I_i}{\ell_i} & 0 & -\frac{6 E_i I_i}{\ell_i^2} & \frac{2 E_i I_i}{\ell_i}\\
	
	& \multicolumn{2}{c}{\multirow{2}{*}{\rotatebox{-30}{sym.}}} & \frac{E_i a_i}{\ell_i} & 0 & 0\\
	
	&  &  &  & \frac{12 E_i I_i}{\ell_i^3} & -\frac{6 E_i I_i}{\ell_i^2}\\
	
	&  &  &  &  & \frac{4 E_i I_i}{\ell_i}
	\end{pmatrix}
	\end{equation}
\end{figure*}
\addtocounter{equation}{-3}

Let us consider the problem of optimum design of frame structures composed a finite number of nodes, $n_\mathrm{n}$, and of admissible elements, $n_\mathrm{e}$, defining the so-called ground structure~\cite{Dorn1964}. The frame elements are attributed with non-negative, and thus possibly zero, given-shaped cross-sectional areas $\mathbf{a}$. For simplicity, we assume, in the following text, that the moments of inertia of individual cross-sections $I_i$ are polynomials of degree two\footnote{The same procedure can be employed for higher-degree polynomials, so that all cross-sectional parameters can be optimized concurrently. We restrict ourselves, however, to the polynomials of degree two to maintain a~simpler notation.}, i.e., $I_i = c_{\mathrm{I},i} a_i^2$ for some $c_{\mathrm{I},i} > 0$, which includes all cross-sections with given aspect ratios of all their components. Note that in such a~case, each cross-section is fully determined by its area. Optimizing the cross-sectional areas $\mathbf{a}$ we search the maximum stiff structures within the available volume $\overline{V}$ of a~linear-elastic material. 

The structural stiffness is measured (inversely) by the compliance $c$, work done by external forces, $\mathbf{f}(\mathbf{a})$,
\begin{equation}\label{eq:comp}
c(\mathbf{a}) = \mathbf{f} (\mathbf{a})^\mathrm{T} \mathbf{u} = \mathbf{u}^\mathrm{T} \mathbf{K}(\mathbf{a}) \mathbf{u},
\end{equation}
where $\mathbf{u}$ constitutes the generalized displacement vector, and $\mathbf{K}(\mathbf{a})$ is the symmetric positive semi-definite stiffness matrix---
a~polynomial function of $\mathbf{a}$,
\begin{equation}
\mathbf{K}(\mathbf{a}) = \sum_{i=1}^{n_\mathrm{e}} \mathbf{K}_i (a_i),
\end{equation}%
\addtocounter{equation}{1}%
assembled from contributions of individual elements $\mathbf{K}_i (a_i)$.
The lower the compliance, the stiffer is the structure with respect to the external forces. 

In this contribution, we use the element stiffness matrix of Euler--Bernoulli frame elements \eqref{eq:stiffmat}, with $E_i$ denoting the Young modulus. In \eqref{eq:comp}, $\mathbf{f}(\mathbf{a})$ is the external force column vector---a~linear function of $\mathbf{a}$---to allow for self-weight, assembled from the contributions of elements $\mathbf{f}(a_i)$,
\begin{equation}
\mathbf{f}(\mathbf{a}) = \sum_{i=1}^{n_\mathrm{e}} \mathbf{f}_i (a_i).
\end{equation}
In the following text, we assume that $\forall \mathbf{a}>\mathbf{0}: \mathbf{K}(\mathbf{a}) \succ \mathbf{0}$, i.e., the structure is not a kinematic mechanism and the stiffness matrix is positive definite (which is denoted by ``$\succ \mathbf{0}$'') for all positive cross-sectional areas. Since $\mathbf{K}(\mathbf{a})$ has therefore the full rank, we also have $\mathbf{f}(\mathbf{a}) \in \mathrm{Im}\left(\mathbf{K}(\mathbf{a})\right)$, where $\mathrm{Im}(\bullet)$ is the image of $\bullet$.

This optimization problem is formalized as
\begin{subequations}\label{eq:original}
\begin{alignat}{4}
	\min_{\mathbf{a}, \mathbf{u}}\;\; && \mathbf{f}(\mathbf{a})^\mathrm{T} \mathbf{u}\label{eq:obj}\\
	\mathrm{s.t.}\;\; && \mathbf{K}(\mathbf{a}) \mathbf{u} \;&&=&& \;\mathbf{f}(\mathbf{a}),\label{eq:equilibrium}\\
	&& \bm{\ell}^\mathrm{T} \mathbf{a} \;&&\le&&\; \overline{V},\label{eq:vol}\\
	&& \mathbf{a} \;&&\ge&&\; \mathbf{0},\label{eq:neg}
\end{alignat}
\end{subequations}
with $\bm{\ell}$ being the column vector of the frame elements lengths. Notice that in general, \eqref{eq:original} constitutes a~non-convex non-linear optimization problem because of the bilinear objective function~\eqref{eq:obj} and the polynomial equilibrium equality~\eqref{eq:equilibrium} with possibly singular $\mathbf{K}(\mathbf{a})$. On the other hand, the volume constraint~\eqref{eq:vol} and the cross-sectional areas non-negativity constraint~\eqref{eq:neg} are affine functions of the design variables, $\mathbf{u}$ and $\mathbf{a}$, and are in turn convex. The optimization problem \eqref{eq:original} can be solved to local optimality using standard numerical optimization techniques. In particular, we will solve this problem using the interior-point method implemented in the \texttt{fmincon} function of \textsc{Matlab}.

\subsection{Optimality criteria}

The inherent difficulty of singularity of $\mathbf{K}(\mathbf{a})$ in \eqref{eq:equilibrium} can be circumvented by assuming a~small positive lower-bound on the cross-sectional areas~\cite{Bendsoe2003}, which is denoted by $\varepsilon$ in this study. Consequently, the former topology optimization problem \eqref{eq:original} is effectively transformed into the sizing one and the displacement field $\mathbf{u}$ can be excluded from the design variables\footnote{At the price of solving the equilibrium equation in each iteration.}. Notice that in the limit when $\varepsilon \rightarrow 0$, these problems are equivalent, but the smaller $\varepsilon$ the higher the condition number of $\mathbf{K}(\mathbf{a})$, and thus it is more difficult to solve \eqref{eq:equilibrium}. Hence, we assume $\varepsilon=10^{-6}$ in this study.

Considering \eqref{eq:original} with $\varepsilon\mathbf{1} \le\mathbf{a}$, its Lagrangian function reads as
\begin{equation}\label{eq:lagrangian}
\begin{split}
\mathcal{L}(\mathbf{a}, \mathbf{u}, \bm{\lambda}, \mu, \bm{\nu}) = \mathbf{f}(\mathbf{a})^\mathrm{T}\mathbf{u} + \bm{\lambda}^\mathrm{T} \left[\mathbf{f}(\mathbf{a}) - \mathbf{K}(\mathbf{a}) \mathbf{u}\right]\\
+ \mu \left( \bm{\ell}^\mathrm{T} \mathbf{a} - \overline{V}\right) + \bm{\nu}^\mathrm{T}\left(\varepsilon\mathbf{1}-\mathbf{a}\right),
\end{split}
\end{equation}
with the Lagrange multipliers $\bm{\lambda}$, $\mu$, and $\bm{\nu}$; and $\mathbf{1}$ denoting the vector of all ones. In addition to the primal feasibility \eqref{eq:equilibrium}--\eqref{eq:neg}, the Karush--Kuhn--Tucker conditions require feasibility of the dual and complementary slackness,
\begin{subequations}\label{eq:lmult}
	\begin{align}
	\mu &\ge 0,\\
	\mu\left( \bm{\ell}^\mathrm{T} \mathbf{a} -\overline{V}\right) &= 0,\\
	\bm{\nu} &\ge \mathbf{0},\label{eq:lmult3}\\
	\bm{\nu}^\mathrm{T}\left(\varepsilon\mathbf{1}-\mathbf{a}\right) &= \mathbf{0},\label{eq:lmult4}
	\end{align}
\end{subequations}
together with the stationarity of the Lagrangian with respect to $\mathbf{u}$,
\begin{equation}\label{eq:lammult}
\mathbf{0} = \diff{\mathcal{L}(\mathbf{a}, \mathbf{u}, \bm{\lambda}, \mu, \bm{\nu})}{\mathbf{u}} = \left[ \mathbf{f}(\mathbf{a}) - \mathbf{K}(\mathbf{a})\bm{\lambda} \right]^\mathrm{T},
\end{equation}
and $\mathbf{a}$. However, because the equation $\mathbf{K}(\mathbf{a}) \bm{\lambda} = \mathbf{f}(\mathbf{a})$ in \eqref{eq:lammult} possesses a unique solution due to $\mathbf{K}(\mathbf{a}) \succ \mathbf{0}$, we have $\bm{\lambda} = \mathbf{u}$. Using this observation, the necessary first-order optimality conditions read as
\begin{equation}
0 = \diffp{\mathcal{L}}{{a_i}} = 2\mathbf{u}^\mathrm{T}\diffp{\mathbf{f}(\mathbf{a})}{{a_i}} - \mathbf{u}^\mathrm{T} \diffp{\mathbf{K}(\mathbf{a})}{{a_i}} \mathbf{u} + \mu \ell_i - \nu_i.\label{eq:optcond}
\end{equation}
Conditions \eqref{eq:optcond} and \eqref{eq:lmult4} then imply that at the optimum, the frame elements with $a_i > \varepsilon$ have equal constant energy
\begin{equation}\label{eq:mutual}
\mu = \frac{1}{\ell_i} \mathbf{u}^\mathrm{T} \diffp{\mathbf{K}(\mathbf{a})}{{a_i}} {\mathbf{u}} 
-\frac{2}{\ell_i}\mathbf{u}^\mathrm{T}\diffp{\mathbf{f}(\mathbf{a})}{{a_i}}.
\end{equation}

Aiming to satisfy \eqref{eq:mutual}, optimality criteria methods~\cite{Bendsoe2003} build update schemes which increase stiffnesses of elements with energies higher than $\mu$, and, conversely, decrease stiffnesses of elements with the energy lower than $\mu$. The levels are balanced in an iterative process, based on the value of $\mu$. In each iteration, the relative change of the design variables is bounded by the move limit $\zeta$, assumed as $\zeta= 0.2$ in this paper. Consequently, the fix-point update scheme reads as
\begin{equation}\label{eq:ocupdate}
\begin{split}
a_i^{(k+1)} = \max\left\{\max\left\{(1-\zeta)a_i^{(k)},\varepsilon\right\},\right.\\
\left.a_i^{(k)} \left[ b_i^{(k)} \right]^\eta
\right\}
\end{split}
\end{equation}
with the tuning parameter $\eta = 0.3$ and with
\begin{equation}\label{eq:ocb}
b^{(k)}_i = \frac{\mathbf{u}^\mathrm{T} \diffp{\mathbf{K}(\mathbf{a})}{{a_i}} \mathbf{u} - 2 \mathbf{u}^\mathrm{T} \diffp{\mathbf{f}(\mathbf{a})}{{a_i}}
}{\mu\ell_i}.
\end{equation}
Clearly, if $\forall i \in \{1,\;\dots\;,n_\mathrm{e}\}: b^{(k)}_i = 1$, we reach a~local minimum as \eqref{eq:mutual} is satisfied.

Combination of \eqref{eq:ocupdate}, \eqref{eq:ocb} with \eqref{eq:vol} allows us to write the (current) volume $V = \bm{\ell}^\mathrm{T} \mathbf{a}^{(k+1)}$ as a continuous function of the multiplier $\mu$. It can be seen from \eqref{eq:ocb} that $V$ is a non-increasing function of $\mu$. In fact, strict decreasing occurs when $\varepsilon \mathbf{I} < \mathbf{a}$. Consequently, the bisection algorithm is used to find $\mu$ such that the volume constraint is satisfied.

\subsection{Nonlinear semidefinite programming}\label{sec:nsdp}

In this section, we describe another approach to eliminate the displacement field variables $\mathbf{u}$ from the optimization problem formulation. First, let us rewrite \eqref{eq:original} as
\begin{subequations}\label{eq:sdp_prior}
	\begin{alignat}{4}
	\min_{\mathbf{a},c}\;\; && c\qquad\qquad\qquad\qquad\;\,\label{eq:nsdpO}\\
	\mathrm{s.t.}\;\; && c - \mathbf{f}(\mathbf{a})^\mathrm{T} \mathbf{K}(\mathbf{a})^{\dagger} \mathbf{f}(\mathbf{a})  \;&&=&&\; 0,\label{eq:sdporig}\\
	&& \bm{\ell}^\mathrm{T} \mathbf{a} \;&&\le&&\; \overline{V},\\
	&& \mathbf{a}\;&&\ge&&\; \mathbf{0},
	\end{alignat}
\end{subequations}
where $\mathbf{K}(\mathbf{a})^\dagger$ denotes the Moore-Penrose pseudo-inverse of $\mathbf{K}(\mathbf{a})$. Because we require that $\forall \mathbf{a} > \mathbf{0}:\mathbf{K}(\mathbf{a}) \succ \mathbf{0}$, the (possible) singularity of $\mathbf{K}(\mathbf{a})$ is caused exclusively by zero rows and columns belonging to the degrees of freedom without any attached finite element (or, equivalently, $a_i = 0$ for all attached elements in that node). This assumption allows us to partition the stiffness matrix into a positive definite principal submatrix $\hat{\mathbf{K}}(\mathbf{a})$ and zero blocks, so that the pseudo-inverse equals
\begin{equation}
\mathbf{K}(\mathbf{a})^{\dagger} = \begin{pmatrix}
\hat{\mathbf{K}}(\mathbf{a})^{-1} & \mathbf{0}\\
\mathbf{0}^\mathrm{T} & \mathbf{0}
\end{pmatrix}.\label{eq:pinv}
\end{equation}
Using the same partitioning, the force vector $\mathbf{f}(\mathbf{a})$ is split into $\begin{pmatrix}
\hat{\mathbf{f}}(\mathbf{a})^\mathrm{T}& \mathbf{0}^\mathrm{T}
\end{pmatrix}^\mathrm{T},$
in which the term $\mathbf{0}^\mathrm{T}$ appears due to the original assumption that $\mathbf{f}(\mathbf{a}) \in \mathrm{Im}(\mathbf{K}(\mathbf{a}))$. Consequently, we see that \eqref{eq:sdporig} can be rewritten to
\begin{equation}
c - \hat{\mathbf{f}}(\mathbf{a})^\mathrm{T} \hat{\mathbf{K}}(\mathbf{a})^{-1} \hat{\mathbf{f}}(\mathbf{a}) = \mathbf{0},
\end{equation}
and \eqref{eq:sdp_prior} is therefore equivalent to \eqref{eq:original}.

From \eqref{eq:pinv} we have $\mathbf{K}^{\dagger}(\mathbf{a}) \succeq \mathbf{0}$, so that $c\ge 0$ based on \eqref{eq:sdporig}. Moreover, iff $\hat{\mathbf{f}}(\mathbf{a}) \neq \mathbf{0}$, we have both $\mathbf{f}(\mathbf{a})^\mathrm{T} \mathbf{K}(\mathbf{a})^{\dagger} \mathbf{f}(\mathbf{a}) > 0$ and $c>0$. Because we minimize $c$ \eqref{eq:nsdpO} and $\mathbf{f}(\mathbf{a})^\mathrm{T} \mathbf{K}(\mathbf{a})^{\dagger} \mathbf{f}(\mathbf{a})$ is bounded from below, \eqref{eq:sdporig} can be simplified to the one-sided inequality
\begin{equation}
c - \mathbf{f}(\mathbf{a})^\mathrm{T} \mathbf{K}(\mathbf{a})^{\dagger} \mathbf{f}(\mathbf{a}) \ge 0.\label{eq:toschur}
\end{equation}
To use the generalized Schur complement lemma, e.g., \cite[Theorem 16.1]{Gallier2011}, we further need to show that
\begin{equation}
\left[\mathbf{I} - \mathbf{K}(\mathbf{a}) \mathbf{K}(\mathbf{a})^\dagger\right] \mathbf{f}(\mathbf{a}) = \mathbf{0}\label{eq:toschur2}
\end{equation}
holds, with $\mathbf{I}$ denoting the identity matrix. Indeed, \eqref{eq:toschur2} is always satisfied because $\mathbf{f}(\mathbf{a}) \in \mathrm{Im}(\mathbf{K}(\mathbf{a}))$, so that we can substitute $\mathbf{f}(\mathbf{a})$ by $\mathbf{K}(\mathbf{a}) \mathbf{v}$, where $\mathbf{v}$ is a vector of coefficients of the linear combination. Then, \eqref{eq:toschur2} is equivalent to
\begin{equation}
\left[\mathbf{K}(\mathbf{a}) - \mathbf{K}(\mathbf{a}) \mathbf{K}(\mathbf{a})^\dagger \mathbf{K}(\mathbf{a})\right] \mathbf{v} = \mathbf{0},\label{eq:toschur3}
\end{equation}
with the term in the square brackets always zero~\cite[Lemma 14.1]{Gallier2011}, as $\mathbf{K}(\mathbf{a})$ is symmetric.

Finally, application of the generalized Schur complement lemma to \eqref{eq:toschur} and \eqref{eq:toschur3} provides us with
\begin{subequations}\label{eq:nsdp}
	\begin{alignat}{4}
	\min_{\mathbf{a}, c}\;\; && c\qquad\qquad\qquad\qquad\;\;\\
	\mathrm{s.t.}\;\; && 
	\begin{pmatrix}
	c & -\mathbf{f}(\mathbf{a})^\mathrm{T}\\
	-\mathbf{f}(\mathbf{a})^\mathrm{T} & \mathbf{K}(\mathbf{a})
	\end{pmatrix}\;&&\succeq&&\; \mathbf{0},\label{eq:pmi}\\
	&& \bm{\ell}^\mathrm{T} \mathbf{a} \;&&\le&&\; \overline{V},\label{eq:nsdp_vol}\\
	&& \mathbf{a}\;&&\ge&&\; \mathbf{0},\label{eq:nsdp_areas}
	\end{alignat}
\end{subequations}
which is a~non-linear semidefinite program equivalent to both \eqref{eq:sdp_prior} and \eqref{eq:original}. Moreover, if we substitute $c$ with $1/s$, where $0<s<\infty$ is a measure of stiffness, \eqref{eq:nsdp} is, after the application of the (standard) Schur complement lemma, e.g., \cite[Proposition 16.1]{Gallier2011}, reducible to
\begin{subequations}\label{eq:nsdp2}
	\begin{alignat}{4}
	\max_{\mathbf{a}, s}\;\; && s\qquad\qquad\qquad\qquad\;\;\\
	\mathrm{s.t.}\;\; && \mathbf{K}(\mathbf{a}) - s\mathbf{f}(\mathbf{a})\mathbf{f}(\mathbf{a})^\mathrm{T} \;&&\succeq&&\; \mathbf{0},\label{eq:pmi2}\\
	&& \bm{\ell}^\mathrm{T} \mathbf{a} \;&&\le&&\; \overline{V},\\
	&& \mathbf{a}\;&&\ge&&\; \mathbf{0},
	\end{alignat}
\end{subequations}
which is a~useful reformulation when $\mathbf{f}$ is constant, i.e., self-weight is not considered.

It shall be noted that due to the polynomial matrix inequalities \eqref{eq:pmi} and \eqref{eq:pmi2} both the optimization problems are non-convex in general. They can still be solved efficiently (to local optimality) using, e.g., augmented Lagrangian methods~\cite{Kocvara2003,Fiala2013}. In this contribution, we solve these problems using the open-source \textsc{PenLab} optimizer~\cite{Fiala2013}.

\subsection{Polynomial optimization}

Having introduced three different local approaches to frame structure optimization, it is natural and expected that one asks for a global optimization technique. In this section, we exploit the fact that although \eqref{eq:nsdp} is non-convex it is indeed a~polynomial optimization (PO) problem at the same time, which allows us to employ modern PO techniques, namely the Lasserre hierarchy~\cite{Lasserre2001,Lasserre2015}, successively building tighter and tighter convex outer semidefinite programming (SDP) approximations called relaxations~\cite[Corollary 4.3]{Lasserre2015}.

To develop a~more efficient formulation suitable for PO, we first recognize that the design variables in \eqref{eq:nsdp} are all bounded both from below and above:
\begin{subequations}\label{eq:sc}
	\begin{alignat}{6}
	0 &&\le&& a_i &&\le&& \frac{\overline{V}}{\ell_i},&&\quad \forall i\in \{1,\;\dots\;,n_\mathrm{e}\},\label{eq:ac}\\
	0 &&\le&& c\; &&\le&& \hat{c}.\label{eq:cc}&& 
	\end{alignat}
\end{subequations}
While the lower bound in \eqref{eq:ac} is caused by the non-negativity of the cross-sectional areas \eqref{eq:nsdp_areas}, the upper bounds arise from the volume constraint~\eqref{eq:nsdp_vol}: none of the structural elements can occupy larger volume than $\overline{V}$. In the compliance case \eqref{eq:cc}, the lower bound\footnote{In fact the strict inequality $0<c$ is satisfied in all non-trivial optimization problems.} is due to $\mathbf{K}(\mathbf{a}) \succeq \mathbf{0}$, recall the discussion in Section~\ref{sec:nsdp}, and the upper bound is provided by an arbitrary feasible design, e.g., the uniform cross-sectional areas
\begin{equation}
\mathbf{a}=\frac{\overline{V}}{\mathbf{1}^\mathrm{T} \bm{\ell}} \mathbf{1},\label{eq:auni}
\end{equation}
which are used for the upper-bound computation in this paper. Then, the compliance bound is computed from 
\begin{equation}\label{eq:compbound}
\hat{c} = \mathbf{f}(\mathbf{a})^\mathrm{T} \mathbf{K}(\mathbf{a})^{\dagger} \mathbf{f}(\mathbf{a}),
\end{equation}
where $\mathbf{K}(\mathbf{a})^{\dagger}$ reduces to $\mathbf{K}(\mathbf{a})^{-1}$ in the case of \eqref{eq:auni}.

To improve numerical stability of the solution process, we state the optimization problem in terms of the scaled cross-sectional areas $\mathbf{a}_\mathrm{sc}$ and scaled compliance $c_\mathrm{sc}$ instead of the original $\mathbf{a}$ and $c$, where the scaled variables are defined in the $[-1,1]$ domain. Therefore, we have
\begin{equation}
a_i  = \frac{\overline{V} (a_{\mathrm{sc},i} + 1)}{2 \ell_i}
\end{equation}
and
\begin{equation}
c = \frac{\hat{c}  (c_\mathrm{sc}+1)}{2}.
\end{equation}
Rewriting \eqref{eq:sc} in terms of $a_{\mathrm{sc},i}$ and $c_\mathrm{sc}$ as unit ball constraints
\begin{subequations}
	\begin{align}
	a_{\mathrm{sc},i}^2 &\le 1,\quad \forall i \in \{1,\;\dots\;,n_\mathrm{e}\},\\
	c_\mathrm{sc}^2 & \le 1,
	\end{align}
\end{subequations}
we finally arrive at an equivalent formulation to \eqref{eq:nsdp},
\begin{subequations}\label{eq:nsdpSC}
	\begin{alignat}{4}
	\min_{\mathbf{a}_\mathrm{sc}, c_\mathrm{sc}}\;\; && 0.5 \hat{c} (c_\mathrm{sc}+1)\qquad\qquad\qquad\;\;\label{nsdpSC_obj}\\
	\mathrm{s.t.}\;\; && 
	\begin{pmatrix}
	0.5 \hat{c} (c_\mathrm{sc}+1) & -\mathbf{f}(\mathbf{a}_\mathrm{sc})^\mathrm{T}\\
	-\mathbf{f}(\mathbf{a}_\mathrm{sc})^\mathrm{T} & \mathbf{K}(\mathbf{a}_\mathrm{sc})
	\end{pmatrix}\;&&\succeq&&\; \mathbf{0},\label{eq:pmiSC}\\
	&& 2 - n_\mathrm{e} - \mathbf{1}^\mathrm{T} \mathbf{a}_{\mathrm{sc}} \;&&\ge&&\; 0,\label{eq:volSC}\\
	&& 1 - a_{\mathrm{sc},i}^2\;&&\ge&&\; 0,\label{eq:aSC}\\
	&& 1 - c_\mathrm{sc}^2\;&&\ge&&\;0.\label{eq:cSC}
	\end{alignat}
\end{subequations}
In \eqref{eq:nsdpSC}, the objective function \eqref{nsdpSC_obj} and also the constraints \eqref{eq:pmiSC}--\eqref{eq:cSC} are all polynomial inequalities\footnote{The entries in the PMI \eqref{eq:pmiSC} are polynomials, so that the feasible region of the PMI is a semi-algebraic set. Alternatively, we may require the roots of the characteristic polynomial of \eqref{eq:pmiSC} to be non-negative. We refer to \cite{Henrion2006} for more details.} non-negative in the feasible region of \eqref{eq:nsdpSC}, which is therefore a~semi-algebraic set.

\begin{figure*}[!b]
	\begin{subfigure}{0.5\linewidth}
		\begin{tikzpicture}
		\scaling{0.67};
		
		\point{a}{0}{0};
		\point{b}{3}{0};
		\point{d}{7}{0};
		\point{c}{10}{0};
		\point{aD}{0}{-0.50}
		\point{aU}{0}{+0.50}
		\point{bD}{3}{-.50}
		\point{bU}{3}{+0.50}
		\point{bD2}{7}{-.25}
		\point{bU2}{7}{+0.25}
		\point{bt}{3}{+0.40}
		\point{bt2}{3}{-0.40}
		\point{ht}{4.5}{+0.40}
		\point{ht2}{4.5}{-0.40}
		\point{htm}{4.5}{0}
		\point{cD}{10}{-.25}
		\point{cU}{10}{+0.25}
		\point{cH}{11}{0}
		\point{cV}{11}{1}
		\point{dt}{7}{+0.33}
		\point{dt2}{7}{-0.33}
		\point{itm}{5.5}{0}
		\point{it}{5.5}{+0.33}
		\point{it2}{5.5}{-0.33}
		
		\beam{2}{aU}{bU}
		\beam{2}{bU2}{cU}
		\beam{2}{aD}{bD}
		\beam{2}{bD2}{cD}
		\beam{2}{bt}{ht}
		\beam{2}{bt2}{ht2}
		\beam{2}{dt}{it}
		\beam{2}{dt2}{it2}
		
		\beam{3}{a}{htm}
		\beam{3}{itm}{c}
		
		\beam{2}{aD}{aU}
		\beam{2}{bD}{bU}
		\beam{2}{dt2}{dt}
		\beam{2}{cD}{cU}
		
		\support{3}{a}[270]
		
		\dimensioning{1}{a}{b}{-1.6}[$\frac{1}{n_\mathrm{e}}$]
		\dimensioning{1}{d}{c}{-1.6}[$\frac{1}{n_\mathrm{e}}$]
		\dimensioning{1}{b}{d}{-1.6}[$1-\frac{2}{n_\mathrm{e}}$]
		\dimensioning{1}{a}{c}{-2.2}[$1$]
		
		\notation{1}{a}{\circled{$1$}}[below=6mm];
		\notation{1}{b}{\circled{$2$}}[below=6mm];
		\notation{1}{d}{\circled{$n_\mathrm{e}$}}[below=6mm];
		\notation{1}{c}{\circled{$n_\mathrm{n}$}}[below=6mm];
		\notation{4}{a}{b}[$1$];
		\notation{4}{d}{c}[$n_\mathrm{e}$];
		
		\load{1}{c}[45][1.0][0.0]
		\notation{1}{c}{$1$}[above right=9mm];
		\addon{3}{c}{cH}{cV}[-1];
		\notation{1}{c}{$30^\circ$}[right=4mm];
		
		\draw[->] (0.0,0)--(1,0);
		\node at (0.6, -0.15) {$x$};
		\draw[->] (0.0,0)--(0,1);
		\node at (0.15, 0.6) {$y$};
		
		\end{tikzpicture}
		\caption{}
		\label{beamA}
	\end{subfigure}%
	\begin{minipage}{0.5\linewidth}
		\begin{subfigure}{0.3\linewidth}
			\includegraphics[width=\linewidth]{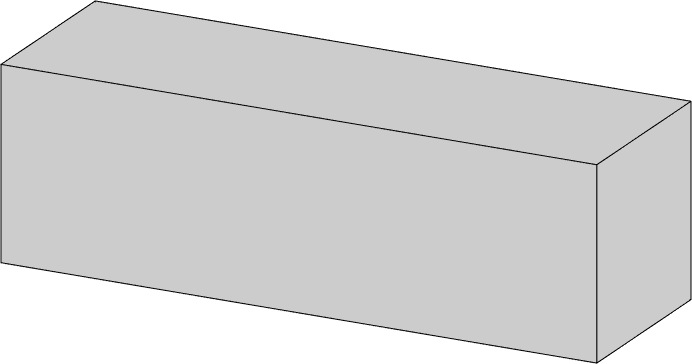}
			\caption{$c^* = 107.50$}
			\label{beamB}
		\end{subfigure}%
		\hfill\begin{subfigure}{0.3\linewidth}
			\includegraphics[width=\linewidth]{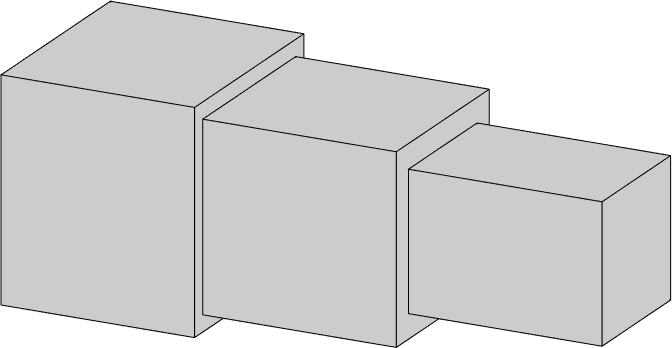}
			\caption{$c^* = 80.30$}
			\label{beamC}
		\end{subfigure}%
		\hfill\begin{subfigure}{0.3\linewidth}
			\includegraphics[width=\linewidth]{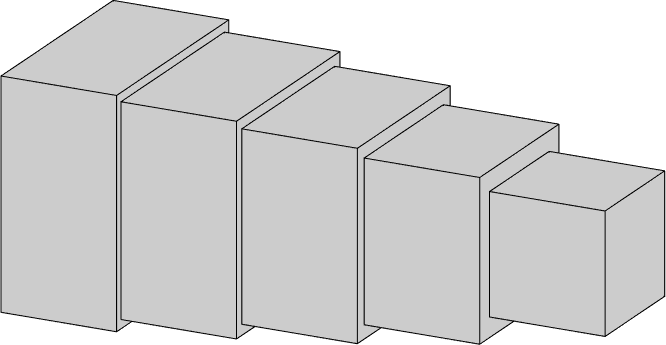}
			\caption{$c^* = 77.19$}
			\label{beamD}
		\end{subfigure}\\
		\begin{subfigure}{0.3\linewidth}
			\includegraphics[width=\linewidth]{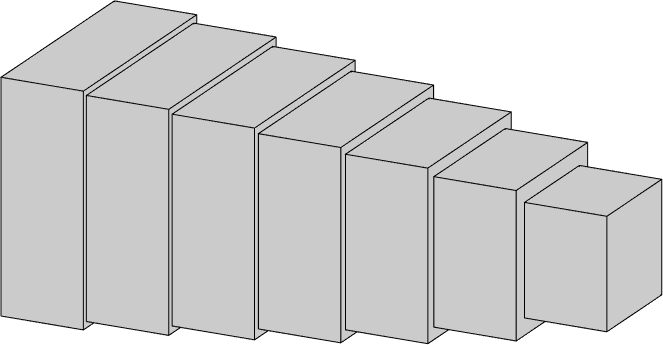}
			\caption{$c^* = 76.23$}
			\label{beamE}
		\end{subfigure}%
		\hfill\begin{subfigure}{0.3\linewidth}
			\includegraphics[width=\linewidth]{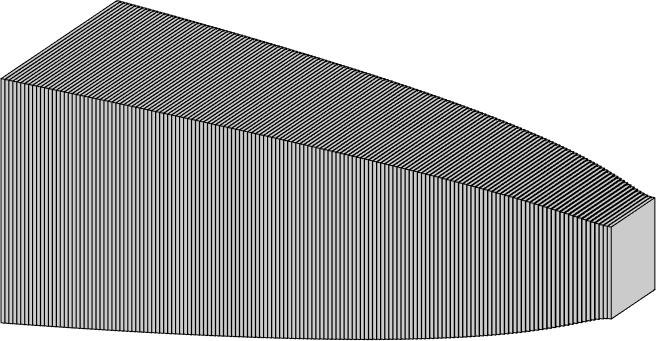}
			\caption{$c=75.19$}
			\label{beamF}
		\end{subfigure}%
		\hfill\begin{subfigure}{0.3\linewidth}
			\includegraphics[width=\linewidth]{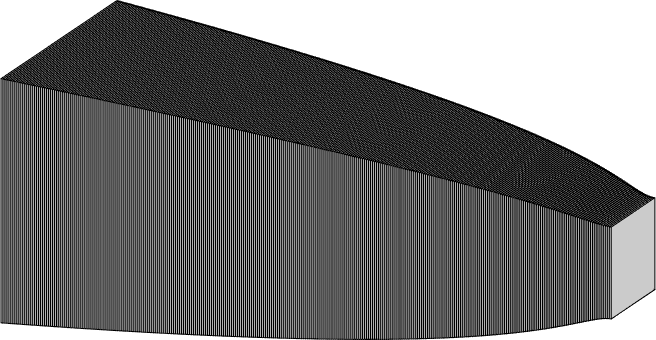}
			\caption{$c=75.19$}
			\label{beamG}
		\end{subfigure}
	\end{minipage}
	\begin{subfigure}{0.33\linewidth}
		\begin{tikzpicture}
		\begin{axis}[
		xlabel={Relaxation order $r$},
		ylabel={Compliance $c$},
		xmin=1, xmax=2.1,
		ymin=0, ymax=140,
		xtick={1,2},
		ytick={0,20,40,60,80,100,120},
		legend pos=south east,
		axis lines=left,
		width=\linewidth,
		]
		\node [above right] at (rel axis cs:0.0,0.65) {\includegraphics[width=1.2cm]{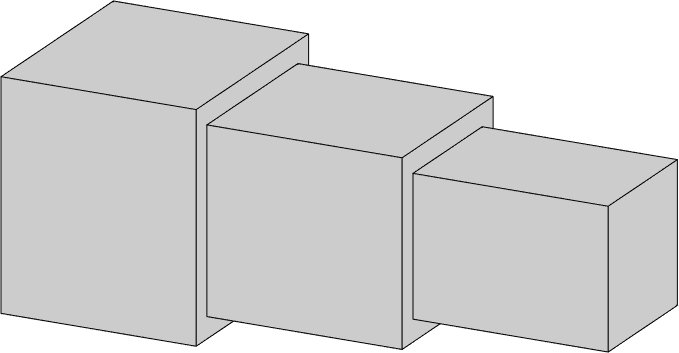}};
		\node [above right] at (rel axis cs:0.65,0.65) {\includegraphics[width=1.2cm]{./include/cantilever/3.png}};
		\addplot[
		color=black, 
		style=dashed,
		mark=triangle,
		mark options={solid},
		] coordinates {
			(1, 80.72)
			(2, 80.30)
		};
		\addlegendentry{$\hat{c}_r$}
		\addplot[
		color=black,
		mark=square,
		mark options={solid},
		] coordinates {
			(1, 35.81)
			(2, 80.30)
		};
		\addlegendentry{$\tilde{c}_r$}
		\addplot[%
		mark=text,
		text mark={},
		text mark as node=true,
		text mark style={star, star point height=2mm,color=blue,fill=white,scale=0.4, draw}
		] coordinates {
			(2, 80.30)
		};
		\addlegendentry{$c^*$}
		\end{axis}
		\end{tikzpicture}
		\caption{}
		\label{beamH}
	\end{subfigure}%
	\begin{subfigure}{0.33\linewidth}
		\begin{tikzpicture}
		\begin{axis}[
		xlabel={Relaxation order $r$},
		ylabel={Compliance $c$},
		xmin=1, xmax=3.1,
		ymin=0, ymax=140,
		xtick={1,2,3},
		ytick={0,20,40,60,80,100,120},
		legend pos=south east,
		axis lines=left,
		width=\linewidth,
		]
		\node [above right] at (rel axis cs:0.0,0.65) {\includegraphics[width=1.2cm]{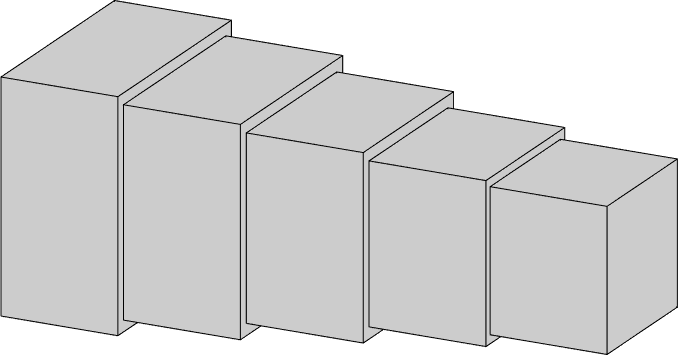}};
		\node [above right] at (rel axis cs:0.325,0.65) {\includegraphics[width=1.2cm]{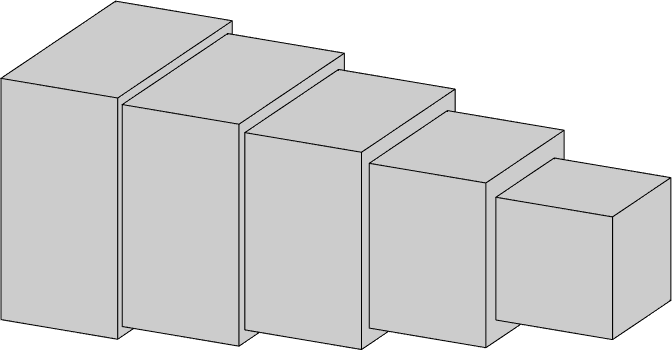}};
		\node [above right] at (rel axis cs:0.65,0.65) {\includegraphics[width=1.2cm]{./include/cantilever/5.png}};
		\addplot[
		color=black, 
		style=dashed,
		mark=triangle,
		mark options={solid},
		] coordinates {
			(1, 79.09)
			(2, 77.37)
			(3, 77.19)
		};
		\addlegendentry{$\hat{c}_r$}
		\addplot[
		color=black,
		mark=square,
		mark options={solid},
		] coordinates {
			(1, 24.72)
			(2, 76.34)
			(3, 77.19)
		};
		\addlegendentry{$\tilde{c}_r$}
		\addplot[%
		mark=text,
		text mark={},
		text mark as node=true,
		text mark style={star, star point height=2mm,color=blue,fill=white,scale=0.4, draw}
		] coordinates {
			(3, 77.19)
		};
		\addlegendentry{$c^*$}
		\end{axis}
		\end{tikzpicture}
		\caption{}
		\label{beamI}
	\end{subfigure}%
	\begin{subfigure}{0.33\linewidth}
		\begin{tikzpicture}
		\begin{axis}[
		xlabel={Relaxation order $r$},
		ylabel={Compliance $c$},
		xmin=1, xmax=3.1,
		ymin=0, ymax=140,
		xtick={1,2,3},
		ytick={0,20,40,60,80,100,120},
		legend pos=south east,
		axis lines=left,
		width=\linewidth,
		]
		\node [above right] at (rel axis cs:0.0,0.65) {\includegraphics[width=1.2cm]{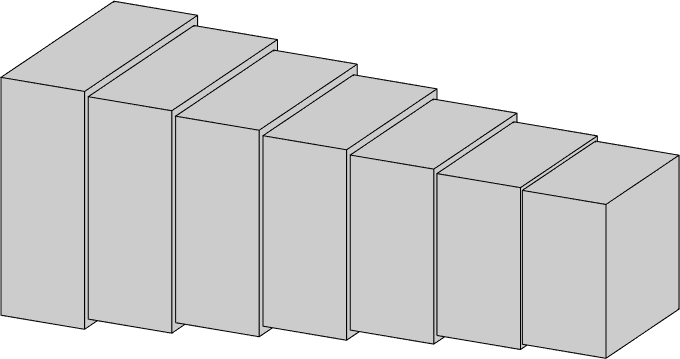}};
		\node [above right] at (rel axis cs:0.325,0.65) {\includegraphics[width=1.2cm]{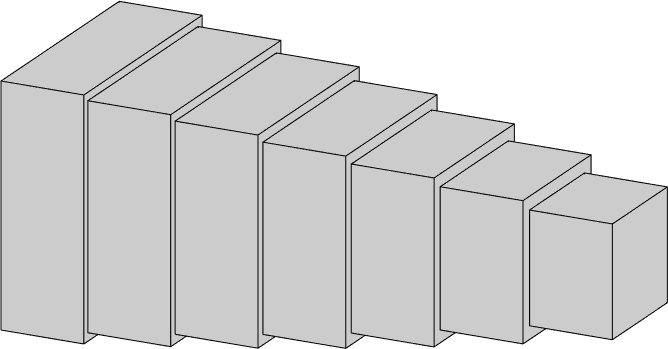}};
		\node [above right] at (rel axis cs:0.65,0.65) {\includegraphics[width=1.2cm]{./include/cantilever/7.png}};
		\addplot[
		color=black, 
		style=dashed,
		mark=triangle,
		mark options={solid},
		] coordinates {
			(1, 79.81)
			(2, 77.02)
			(3, 76.23)
		};
		\addlegendentry{$\hat{c}_r$}
		\addplot[
		color=black,
		mark=square,
		mark options={solid},
		] coordinates {
			(1, 20.00)
			(2, 71.69)
			(3, 76.23)
		};
		\addlegendentry{$\tilde{c}_r$}
		\addplot[%
		mark=text,
		text mark={},
		text mark as node=true,
		text mark style={star, star point height=2mm,color=blue,fill=white,scale=0.4, draw}
		] coordinates {
			(3, 76.23)
		};
		\addlegendentry{$c^*$}
		\end{axis}
		\end{tikzpicture}
		\caption{}
		\label{beamJ}
	\end{subfigure}
	\caption{Cantilever beam design problem. (a) shows the design domain, (b)--(e) are the optimal topologies for $1$, $3$, $5$, and $7$ elements; (f) and (g) display optimized topologies computed by OC with discretization by $150$ and $300$ elements. Figures (h)--(j) show the convergance of PO for $3$, $5$, and $7$ elements.}
\end{figure*}

The optimization problem \eqref{eq:nsdpSC} can readily be solved by building the (Lasserre) hierarchy of convex linear SDP relaxations, with a monotonously converging objective functions to the global optimum. The hierarchy is generated by the \textsc{GloptiPoly}~\cite{Henrion2009} software package interfaced with the \textsc{Yalmip}~\cite{Lofberg2004} toolbox, and the underlying SDP relaxations are solved by the state-of-the-art \textsc{Mosek} optimizer~\cite{mosek}.

\subsubsection{Solution process}

To simplify the notation, let us now denote the objective function polynomial \eqref{nsdpSC_obj} by $p_0$, and the constraining polynomials \eqref{eq:volSC} to \eqref{eq:cSC} by $p_1$ to $p_3$, respectively. Further, let $P_4$ denote the PMI \eqref{eq:pmiSC} and let
\begin{equation}
\mathbf{x} = \begin{pmatrix}c_\mathrm{sc} & a_{\mathrm{sc},1} & \dots & a_{\mathrm{sc},n_\mathrm{e}} \end{pmatrix}^\mathrm{T}
\end{equation}
be the vector of the design variables. Moreover, let
\begin{equation}
\begin{split}
\mathbf{b}_r (\mathbf{x}) = \left(1\;\; x_1\;\;
x_2\;\;\dots\;\;x_{n_\mathrm{e}+1}\;\;
x_1^2 \;\; x_1 x_2\;\;\dots\;\; \right.\\ 
\left. x_1 x_{n_\mathrm{e}+1} \;\; x_2^2 \;\; x_2 x_3\;\; \dots \;\; x_{n_\mathrm{e}+1}^2
\dots \right.\\ 
\left.x_1^r \;\; \dots \;\;x_{n_\mathrm{e}+1}^r\right)^\mathrm{T}
\end{split}
\end{equation}
denote the polynomial space basis of the maximum degree $r$. Then, we can express each of the polynomials $p_j, j \in \{0,\;\dots\;,3\}$ as a linear combination 
\begin{equation}\label{eq:lincombmon}
p_j (\mathbf{x}) = \sum_{\beta=1}^{\lvert\mathbf{b}_r(\mathbf{x}) \rvert} p_{\alpha, j,\beta} (x^\alpha)_\beta
\end{equation}
of monomials 
\begin{equation}
x^\alpha = \prod_{m=1}^{n_\mathrm{e}+1} x_m^{\alpha_m}, \quad \sum_{m=1}^{n_\mathrm{e}+1} \alpha_m \le r,
\end{equation}
indexed in the basis $\mathbf{b}_r(\mathbf{x})$. The vectors $\bm{\alpha}$ associate a~non-negative integer number with each element in $\mathbf{x}$, and the vector $\mathbf{p}_{\alpha, j}$ contains coefficients of the linear combinations of the monomials. Notice that $\sum_{m=1}^{n_\mathrm{e}+1} \alpha_m \le d_j$, where $d_j$ stands for the degree of the polynomial $p_j$. A similar approach is also applied to the elements of $P_4$~\cite{Henrion2006}.

\begin{table*}[!t]
	\begin{subtable}{0.4\linewidth}
		\begin{tabular}{lcccc}
			& $t$ [s] & $\hat{c}$ & $\tilde{c}$ & $a_1$\\
			\hline
			\texttt{fmincon} & $0.7$ & $107.50$ & - & $0.100$ \\
			OC & $0.1$ & $107.50$ & - & $0.100$\\
			NSDP & $0.2$ & $107.50$ & - & $0.100$ \\
			PO$^{(1)}$ & $0.3$ & $\textbf{107.50}$ & $107.50$ & $0.100$
		\end{tabular}
		\caption{$1$ element}
	\end{subtable}%
	\hfill\begin{subtable}{0.52\linewidth}
		\begin{tabular}{lcccccc}
			& $t$ [s] & $\hat{c}$ & $\tilde{c}$ & $a_1$ & $a_2$ & $a_3$\\
			\hline
			\texttt{fmincon} & $0.6$ & $80.30$ & - & $0.142$ & $0.102$ & $0.056$ \\
			OC & $0.1$ & $80.30$ & - & $0.142$ & $0.102$ & $0.056$ \\
			NSDP & $0.3$ & $80.30$ & - & $0.142$ & $0.102$ & $0.056$ \\
			PO$^{(1)}$ & $0.3$ & $80.72$ & $35.81$ & $0.147$ & $0.097$ & $0.056$ \\
			PO$^{(2)}$ & $0.3$ & $\textbf{80.30}$ & $80.30$ & $0.142$ & $0.102$ & $0.056$
		\end{tabular}
		\caption{$3$ elements}
	\end{subtable}\\
	\begin{subtable}{\linewidth}
		\centering
		\begin{tabular}{lcccccccc}
			& $t$ [s] & $\hat{c}$ & $\tilde{c}$ & $a_1$& $a_2$ & $a_3$ & $a_4$ & $a_5$\\
			\hline
			\texttt{fmincon} & $0.4$ & $77.19$ & - & $0.151$ & $0.128$ & $0.103$ & $0.075$ & $0.043$ \\
			OC & $0.1$ & $77.19$ & - & $0.151$ & $0.128$ & $0.103$ & $0.075$ & $0.043$ \\
			NSDP & $0.5$ & $77.19$ & - & $0.151$ & $0.128$ & $0.103$ & $0.075$ & $0.043$ \\
			PO$^{(1)}$ & $0.3$ & $79.09$ & $24.72$ & $0.151$ & $0.123$ & $0.096$ & $0.073$ & $0.058$ \\
			PO$^{(2)}$ & $0.7$ & $77.37$ & $76.34$ & $0.154$ & $0.131$ & $0.103$ & $0.072$ & $0.040$ \\
			PO$^{(3)}$ & $26.8$ & $\textbf{77.19}$ & $77.19$ & $0.151$ & $0.128$ & $0.103$ & $0.075$ & $0.043$
		\end{tabular}
		\caption{$5$ elements}
	\end{subtable}
	\begin{subtable}{\linewidth}
		\centering
		\begin{tabular}{lcccccccccc}
			& $t$ [s] & $\hat{c}$ & $\tilde{c}$ & $a_1$ & $a_2$ & $a_3$ & $a_4$ & $a_5$ & $a_6$ & $a_7$\\
			\hline
			\texttt{fmincon} & $0.5$ & $76.23$ & - & $0.155$ & $0.139$ & $0.122$ & $0.104$ & $0.084$ & $0.061$ & $0.036$ \\
			OC & $0.1$ & $76.23$ & - & $0.155$ & $0.139$ & $0.122$ & $0.104$ & $0.084$ & $0.061$ & $0.036$ \\
			NSDP & $0.4$ & $76.23$ & - & $0.155$ & $0.139$ & $0.122$ & $0.104$ & $0.084$ & $0.061$ & $0.036$ \\
			PO$^{(1)}$ & $0.3$ & $79.81$ & $20.00$ & $0.149$ & $0.130$ & $0.112$ & $0.096$ & $0.081$ & $0.069$ & $0.062$ \\
			PO$^{(2)}$ & $2.5$ & $77.02$ & $71.69$ & $0.166$ & $0.145$ & $0.122$ & $0.099$ & $0.077$ & $0.055$ & $0.036$ \\
			PO$^{(3)}$ & $550.1$ & $\textbf{76.23}$ & $76.23$ & $0.155$ & $0.139$ & $0.122$ & $0.104$ & $0.084$ & $0.061$ & $0.036$
		\end{tabular}
		\caption{$7$ elements}
	\end{subtable}
	\caption{Comparison of the four optimization approaches on the design of the cantilever beam problem. Bold text denotes the proven global optimum.}
	\label{tab:cant}
\end{table*}

\begin{figure*}[!b]
	\begin{subfigure}{0.375\linewidth}
		\centering
		\begin{tikzpicture}
		\scaling{2.3};
		
		\point{a}{0}{0};
		\point{b}{1}{0};
		\point{c}{2}{0};
		\point{d}{0}{1};
		\point{e}{1}{1};
		\point{f}{2}{1};
		
		\beam{2}{a}{b};
		\beam{2}{a}{e};
		\beam{2}{b}{c};
		\beam{2}{b}{d};
		\beam{2}{b}{e};
		\beam{2}{b}{f};
		\beam{2}{c}{e};
		\beam{2}{c}{f};
		\beam{2}{d}{e};
		\beam{2}{e}{f};
		
		\support{3}{a}[270];
		\support{3}{d}[270];
		
		\load{2}{b}[-10][-160][0.6];
		\load{2}{c}[-10][-160][0.6];
		\notation{1}{b}{$1$}[below right=5mm];
		\notation{1}{c}{$2$}[below right=5mm];
		
		\draw[->] (0.0,0)--(1,0);
		\node at (0.6, 0.15) {$x$};
		\draw[->] (0.0,0)--(0,1);
		\node at (0.15, 0.6) {$y$};
		
		\dimensioning{1}{a}{b}{-1.2}[$1$]
		\dimensioning{1}{b}{c}{-1.2}[$1$]
		\dimensioning{2}{a}{d}{-0.8}[$1$]
		
		\notation{4}{a}{b}[$1$];
		\notation{4}{a}{e}[$2$][0.33];
		\notation{4}{b}{c}[$3$];
		\notation{4}{b}{d}[$4$][0.33];
		\notation{4}{b}{e}[$5$];
		\notation{4}{b}{f}[$6$][0.33];
		\notation{4}{c}{e}[$7$][0.33];
		\notation{4}{c}{f}[$8$];
		\notation{4}{d}{e}[$9$];
		\notation{4}{e}{f}[$10$];
		
		\notation{1}{a}{\circled{$1$}}[align=center];
		\notation{1}{b}{\circled{$2$}}[align=center];
		\notation{1}{c}{\circled{$3$}}[align=center];
		\notation{1}{d}{\circled{$4$}}[align=center];
		\notation{1}{e}{\circled{$5$}}[align=center];
		\notation{1}{f}{\circled{$6$}}[align=center];
		\end{tikzpicture}
		\vspace{8mm}
		\caption{}
		\label{fig:10}
	\end{subfigure}%
	\hfill\begin{minipage}{0.175\linewidth}
		\begin{subfigure}{\linewidth}
			\centering
			\includegraphics[width=2.75cm]{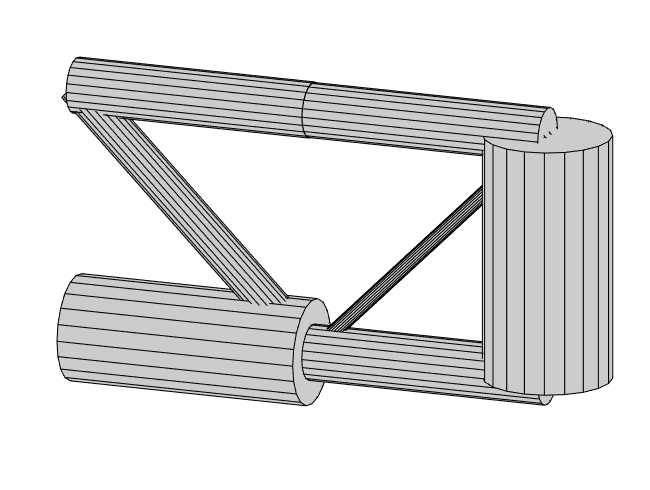}
			\caption{$c = 1042.20$}
			\label{fig:10_lo}
		\end{subfigure}\\
		\vspace{7mm}
		\begin{subfigure}{\linewidth}
			\centering
			\includegraphics[width=2.75cm]{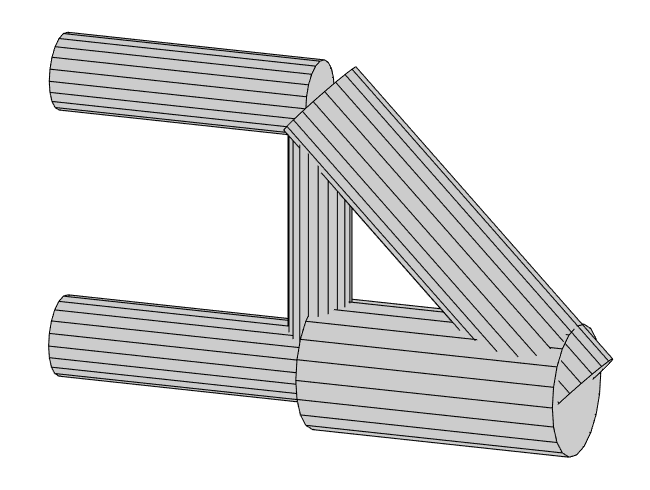}
			\caption{$c^* = 959.32$}
			\label{fig:10_go}
		\end{subfigure}
	\end{minipage}%
	\hfill\begin{subfigure}{0.405\linewidth}
		\begin{tikzpicture}
		\begin{axis}[
		xlabel={Relaxation order $r$},
		ylabel={Compliance $c$},
		xmin=1, xmax=2.1,
		ymin=0, ymax=2350,
		xtick={1,2},
		ytick={0,500,1000,1500,2000},
		legend pos=south east,
		axis lines=left,
		width=\linewidth,
		height=6cm,
		]
		\node [above right] at (rel axis cs:-0.05,0.52) {\includegraphics[width=2.75cm]{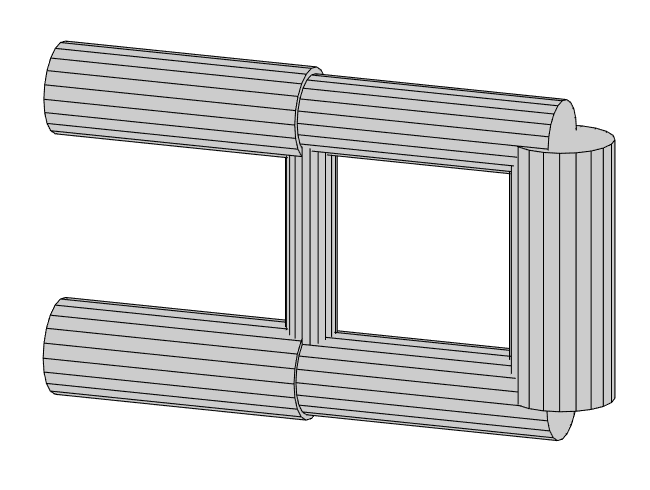}};
		\node [above right] at (rel axis cs:0.45,0.42) {\includegraphics[width=2.75cm]{./include/10frame/psdp10r2_959_3123767.png}};
		\addplot[
		color=black, 
		style=dashed,
		mark=triangle,
		mark options={solid},
		] coordinates {
			(1, 1429.31)
			(2, 959.32)
		};
		\addlegendentry{$\hat{c}_r$}
		\addplot[
		color=black,
		mark=square,
		mark options={solid},
		] coordinates {
			(1, 443.20)
			(2, 959.31)
		};
		\addlegendentry{$\tilde{c}_r$}
		\addplot[%
		mark=text,
		text mark={},
		text mark as node=true,
		text mark style={star, star point height=2mm,color=blue,fill=white,scale=0.4, draw}
		] coordinates {
			(2, 959.315)
		};
		\addlegendentry{$c^*$}
		\end{axis}
		\end{tikzpicture}
		\caption{}
		\label{fig:10_po}
	\end{subfigure}
	\caption{10-beam structure design problem. (a) shows the design domain. While all the tested local optimization algorithms converged to the topology in (b), the global optimum design (c) obtained by PO possesses a clearly different topology. Figure (d) shows the convergence of PO.}
\end{figure*}

Further, introducing $\mathbf{y}$, with its components $y_\beta$ corresponding to a monomial in the basis $\mathbf{b}_r(\mathbf{x})$, we build the Lasserre hierarchy of convex linear semidefinite programming relaxations
\begin{subequations}\label{eq:moment}
	\begin{alignat}{5}
\min_{\mathbf{y}}\;\; && \sum_{\beta=1}^{\lvert \mathbf{b}_r (\mathbf{x}) \rvert} p_{\alpha,0,\beta} y_\beta \;\;\\
\mathrm{s.t.}\;\; && 
\mathbf{M}_r (\mathbf{y})\;&&\succeq&&\; 0,\\
&& \mathbf{M}_{r-d_j} (\mathbf{y}) \;&&\succeq&&\; 0, &&\quad \forall j \in \{1,\;\dots\;,4\},
\end{alignat}
\end{subequations}
and solve it successively with an increasing relaxation order $r$ until the globally optimal solution(s) are found. For all our test cases, this convergence was always finite.


In \eqref{eq:moment}, the matrix $\mathbf{M}_r (\mathbf{y})$ is the moment matrix of the $r$-th order, and $\mathbf{M}_{r-d_j}$ is the $(r-d_j)$-order localization matrix associated with $p_j$ or $P_j$. For more details about these matrices, we refer the reader to \cite{Lasserre2015, Henrion2006}.

\subsubsection{Recognizing global optimality}

There are (at least) two ways to recognize whether $\mathbf{y}_r^*$, a $r$-th order relaxation solution to \eqref{eq:moment}, is globally optimal.

The first (common) way is based a sufficient condition for global optimality of \eqref{eq:moment} \citet{Curto2000,Henrion2005},
\begin{equation}\label{eq:rank}
\mathrm{rank} \;\mathbf{M}_r(\mathbf{y}_r^*) = \mathrm{rank}\; \mathbf{M}_{r-d}(\mathbf{y}^*_r),
\end{equation}
in which $\mathbf{y}^*$ constitutes the optimal solution in the basis $\mathbf{b}_r(c_\mathrm{sc},\mathbf{a}_\mathrm{sc})$, and $\mathrm{rank} \;\mathbf{M}_r(\mathbf{y}^*)$ is less or equal to the number of these solutions~\cite{Lasserre2015}. The values of the original design variables can be extracted using Cholesky or singular value decompositions~\cite{Henrion2005}.

For the optimization problem \eqref{eq:nsdpSC}, we introduce another sufficient condition for global optimality. Let $\tilde{\mathbf{a}}_{r}$  be the (optimized) values of the cross-sectional areas and let $\tilde{c}_{r}$ denote the compliance found in the $r$-th order relaxation. Indeed, $\tilde{c}_{r} \le c^*$ \cite{Lasserre2015}. Moreover, because any feasible design $\mathbf{a}$ generates an upper bound $\hat{c}_{r}$, recall Eq.~\eqref{eq:compbound}, we have
\begin{equation}
\tilde{c}_{r} \le c^* \le \hat{c}_{r}.
\end{equation}
A globally optimal solution to \eqref{eq:nsdpSC} is found if there is no gap,
\begin{equation}\label{eq:necessarry}
\hat{c}_{r} - \tilde{c}_{r} \longrightarrow 0.
\end{equation}
Note that \eqref{eq:necessarry} is substantially simpler to check than \eqref{eq:rank}. If \eqref{eq:necessarry} is not satisfied in the (current) relaxation order $r$, we have at least the quality measure for the current solution in hand.

\begin{table*}[!t]
	\caption{Comparison of the four optimization approaches on the design of $10$-beam structure. Bold text denotes proven global optimum. The cross-sectional area $a_2$ is zero in all test cases.}
	\label{tab:10}
	\begin{tabular}{lcccccccccccc}
		& $t$ [s] & $\hat{c}$ & $\tilde{c}$ & $a_1$ & $a_3$ & $a_4$ & $a_5$ & $a_6$ & $a_7$ & $a_8$ & $a_9$ & $a_{10}$\\
		\hline
		\texttt{fmincon} & $1.0$ & $1042.14$ & - & $0.144$ & $0.040$ & $0.018$ & $0.000$ & $0.003$ & $0.000$ & $0.210$ & $0.039$ & $0.038$ \\
		OC & $0.6$ & $1042.33$ & - & $0.146$ & $0.039$ & $0.017$ & $0.000$ & $0.002$ & $0.000$ & $0.212$ & $0.039$ & $0.037$ \\
		NSDP & $1.0$ & $1042.20$ & - & $0.144$ & $0.040$ & $0.018$ & $0.000$ & $0.003$ & $0.000$ & $0.210$ & $0.039$ & $0.038$ \\
		PO$^{(1)}$ & $0.3$ & $1429.31$ & $443.20$ & $0.103$ & $0.082$ & $0.000$ & $0.029$ & $0.000$ & $0.000$ & $0.121$ & $0.095$ & $0.069$ \\
		PO$^{(2)}$ & $5.2$ & $\textbf{959.32}$ & $959.31$ & $0.070$ & $0.186$ & $0.000$ & $0.043$ & $0.000$ & $0.098$ & $0.000$ & $0.064$ & $0.000$
	\end{tabular}
\end{table*}

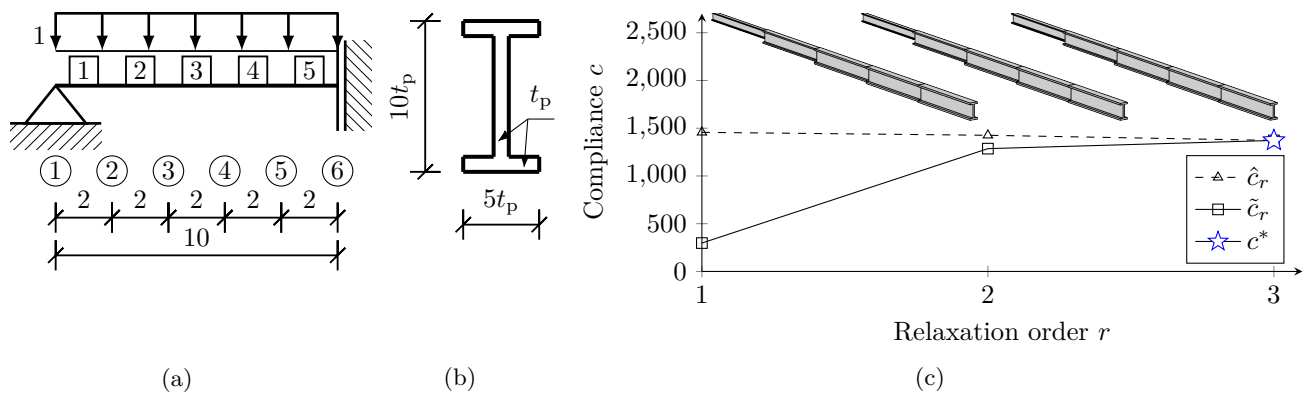
\begin{figure*}[!b]
	\begin{subfigure}{0.27\linewidth}
		\begin{tikzpicture}[baseline=(current bounding box.center)]
		\scaling{0.53};
		
		\point{a}{0}{0};
		\point{aU}{0}{0.3};
		\point{b}{1.4}{0};
		\point{c}{2.8}{0};
		\point{d}{4.2}{0};
		\point{e}{5.6}{0};
		\point{h}{7}{0};
		\point{hU}{7}{0.3};
		\point{lf}{-0.8}{0.8};
		
		\beam{2}{a}{h};
		\notation{4}{a}{b}[$1$];
		\notation{4}{b}{c}[$2$];
		\notation{4}{c}{d}[$3$];
		\notation{4}{d}{e}[$4$];
		\notation{4}{e}{h}[$5$];
		\notation{1}{lf}{$1$};
		
		\lineload{1}{aU}{hU}[0.5][0.5][0.08750];
		
		\support{1}{a};
		\support{4}{h}[90];
		
		\notation{1}{a}{\circled{$1$}}[below=8mm];
		\notation{1}{b}{\circled{$2$}}[below=8mm];
		\notation{1}{c}{\circled{$3$}}[below=8mm];
		\notation{1}{d}{\circled{$4$}}[below=8mm];
		\notation{1}{e}{\circled{$5$}}[below=8mm];
		\notation{1}{h}{\circled{$6$}}[below=8mm];
		
		\dimensioning{1}{a}{b}{-1.75}[$2$];
		\dimensioning{1}{b}{c}{-1.75}[$2$];
		\dimensioning{1}{c}{d}{-1.75}[$2$];
		\dimensioning{1}{d}{e}{-1.75}[$2$];
		\dimensioning{1}{e}{h}{-1.75}[$2$];
		\dimensioning{1}{a}{h}{-2.25}[$10$];
		
		\end{tikzpicture}
		\vspace{10mm}
		\caption{}
		\label{fig:girderA}
	\end{subfigure}%
	\hfill\begin{subfigure}{0.13\linewidth}	
		\begin{tikzpicture}[baseline=(current bounding box.center)]
		\scaling{0.2};
		
		\point{a}{0}{0};
		\point{b}{5}{0};
		\point{c}{5}{1};
		\point{d}{3}{1};
		\point{e}{3}{9};
		\point{f}{5}{9};
		\point{g}{5}{10};
		\point{h}{0}{10};
		\point{i}{0}{9};
		\point{j}{2}{9};
		\point{k}{2}{1};
		\point{l}{0}{1};
		
		\point{t}{2.5}{10};
		
		\point{a1}{4}{0.25};
		\point{a2}{2.25}{2};
		\point{a3}{4.25}{3.5};
		\point{a4}{6.0}{3.5};
		\draw [{Latex}-](a1) -- (a3);
		\draw [{Latex}-](a2) -- (a3);
		\draw (a3) -- node[above]{$t_\mathrm{p}$} (a4);
		
		\beam{2}{a}{b};
		\beam{2}{b}{c};
		\beam{2}{c}{d};
		\beam{2}{d}{e};
		\beam{2}{e}{f};
		\beam{2}{f}{g};
		\beam{2}{g}{h};
		\beam{2}{h}{i};
		\beam{2}{i}{j};
		\beam{2}{j}{k};
		\beam{2}{k}{l};
		\beam{2}{l}{a};

		\dimensioning{1}{a}{b}{-0.7}[$5 t_\mathrm{p}$];
		\dimensioning{2}{b}{g}{-0.5}[$10 t_\mathrm{p}$];
		\end{tikzpicture}
		\vspace{14mm}
		\caption{}
		\label{fig:girderB}
	\end{subfigure}%
	\hfill\begin{subfigure}{0.56\linewidth}
		\begin{tikzpicture}
		\begin{axis}[
		xlabel={Relaxation order $r$},
		ylabel={Compliance $c$},
		xmin=1, xmax=3.1,
		ymin=0, ymax=2700,
		xtick={1,2,3},
		ytick={0,500,1000,1500,2000,2500},
		legend pos=south east,
		axis lines=left,
		width=\linewidth,
		height=5cm,
		]
		\node [above right] at (rel axis cs:0.0,0.55) {\includegraphics[width=3.5cm]{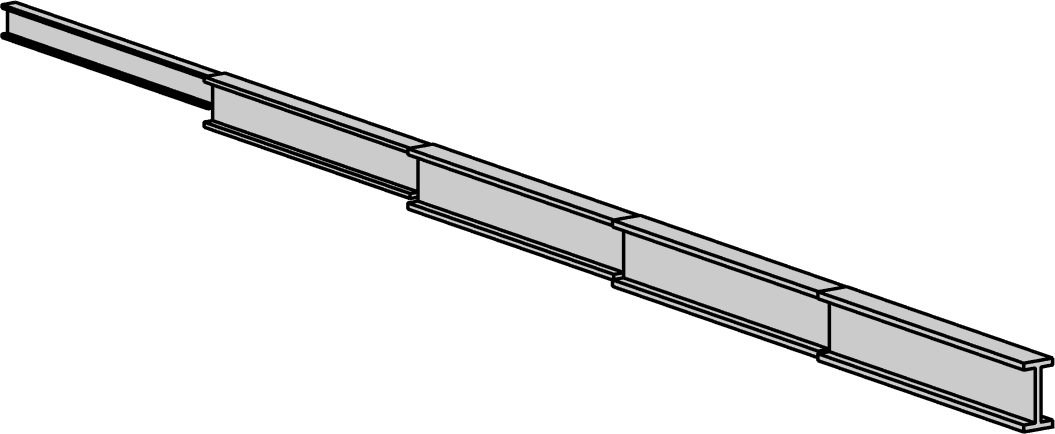}};
		\node [above right] at (rel axis cs:0.25,0.55) {\includegraphics[width=3.5cm]{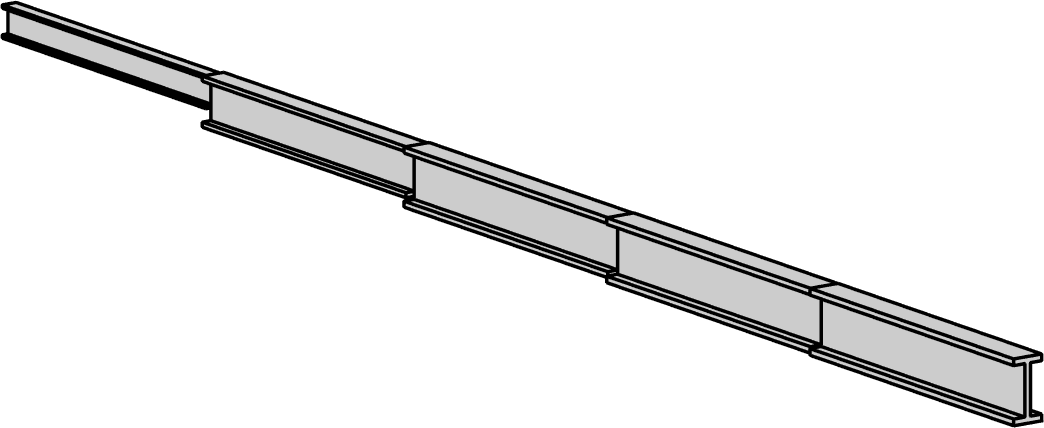}};
		\node [above right] at (rel axis cs:0.5,0.55) {\includegraphics[width=3.5cm]{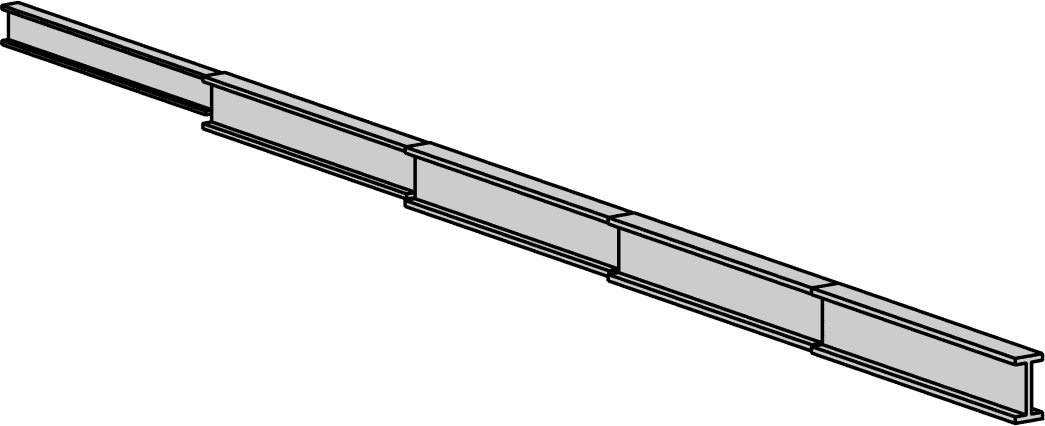}};
		\addplot[
		color=black, 
		style=dashed,
		mark=triangle,
		mark options={solid},
		] coordinates {
			(1, 1456.75)
			(2, 1426.05)
			(3, 1372.25)
		};
		\addlegendentry{$\hat{c}_r$}
		\addplot[
		color=black,
		mark=square,
		mark options={solid},
		] coordinates {
			(1, 297.34)
			(2, 1286.05)
			(3, 1372.25)
		};
		\addlegendentry{$\tilde{c}_r$}
		\addplot[%
		mark=text,
		text mark={},
		text mark as node=true,
		text mark style={star, star point height=2mm,color=blue,fill=white,scale=0.4, draw}
		] coordinates {
			(3, 1372.25)
		};
		\addlegendentry{$c^*$}
		\end{axis}
		\end{tikzpicture}
		\caption{}
		\label{fig:girderC}
	\end{subfigure}
	\caption{Girder beam design problem. (a) shows the design domain, (b) the cross-sectional shape, and (c) the convergence of polynomial optimization.}
	\label{fig:girder}
\end{figure*}

\section{Sample problems}\label{sec:sample}

In this section, we describe and solve three rather small-scaled structural optimization problems, and compare the four optimization approaches.

\subsection{Cantilever beam}

Consider a (simple) cantilever beam design problem, as shown in Fig.~\ref{beamA}. The beam is discretized into $n_\mathrm{e}$ finite elements of equal lengths, each of them assigned a~prismatic square cross-section. The beam is loaded with a skew force at its tip. Further, we assume (dimensionless) $E=1$ and $\overline{V} = 0.1$.

It can be seen from Table~\ref{tab:cant} that the problem is relatively ``easy'' to solve, as all of the tested algorithms converge to the unique global optima (shown in \ref{beamB}--\ref{beamE}), proved by the PO approach (Figs. \ref{beamH}--\ref{beamJ}). In fact, we include this optimization problem mainly to show that PO scales unambiguously worse with the problem dimension than the other methods, and that the optimality criteria (OC) method is the fastest one. Using OC, we can optimize the discretizations of $150$ or $300$ elements in $0.4$ and $0.8$ seconds (Figs.~\ref{beamF} and \ref{beamH}).

\subsection{10-beam frame structure}

\begin{table*}[!t]
	\begin{tabular}{lcccccccc}
		& $t$ [s] & $\hat{c}$ & $\tilde{c}$ & $a_1$ & $a_2$ & $a_3$ & $a_4$ & $a_5$\\
		\hline
		\texttt{fmincon} & $0.4$ & - & - & $0.394$ & $0.256$ & $0.000$ & $0.003$ & $0.000$ \\
		OC & $0.1$ & $1372.25$ & - & $0.010$ & $0.017$ & $0.022$ & $0.025$ & $0.026$\\
		NSDP & $2.0$ & - & - & $0.000$ & $0.000$ & $0.000$ & $0.000$ & $0.000$ \\
		PO$^{(1)}$ & $0.2$ & $1456.75$ & $297.34$ & $0.007$ & $0.015$ & $0.022$ & $0.027$ & $0.029$\\
		PO$^{(2)}$ & $1.0$ & $1426.05$ & $1286.44$ & $0.007$ & $0.016$ & $0.023$ & $0.026$ & $0.028$\\
		PO$^{(3)}$ & $37.7$ & $\textbf{1372.25}$ & $1372.25$ & $0.010$ & $0.017$ & $0.022$ & $0.025$ & $0.026$
	\end{tabular}
	\caption{Comparison of the four optimization approaches on the design of the girder beam problem. Bold text denotes the proven global optimum.}
	\label{tab:girder}
\end{table*}

Second, we consider the topology optimization problem of designing circular cross-sections of the $10$-beam structure shown in Fig.~\ref{fig:10}. This structure is loaded by two moments of magnitudes $1$ and $2$, placed at the nodes \circled{2} and \circled{3}, respectively. Note that there is no intersection between the elements \squared{3} and \squared{4}, and between \squared{6} and \squared{7}. We assume the Young modulus $E=1$ and the volume bound $\overline{V} = 0.5$.

Evaluating the four optimization methods, Table~\ref{tab:10}, we see that while the local optimization methods converged to nearby local minima of topologies shown in Fig.~\ref{fig:10_lo}, PO found the global optimum of a clearly different topology, Fig.~\ref{fig:10_go}. Moreover, the global optimum possesses the compliance lower by $8\%$, which was reached in the second relaxation, Fig.~\ref{fig:10_po}.

\subsection{I-shaped girder with self-weight}

Last, we consider a simply-supported I-shaped girder beam of the span $20$, loaded by a uniform load $1$ and self-weight. Due to the problem symmetry, we consider only one half of the problem, Fig.~\ref{fig:girderA}, and discretize it using $5$ finite elements of equal length. This (half of the) girder beam is allowed to utilize at most $\overline{V} = 0.2$ material of the Young modulus $E=10^4$ and of density $\rho = 3$. The beam cross-sections are parameterized by the parameter $t_\mathrm{p}$, which denotes the thicknesses of the flanges and web. The cross-sectional area equals $a(t_\mathrm{p}) = 18 t_\mathrm{p}^2$. We assume that the upper surfaces of the top flanges are aligned, so that we have $I (t_\mathrm{p}) = 696 t_\mathrm{p}^4 = 58/27 a(t_\mathrm{p})^2$. Consequently, we optimize $a_i$ of individual elements, rather than $t_\mathrm{p}$.

Table \ref{tab:girder} reveals that for this specific problem the OC converged to the global optimum. On the other hand, the \texttt{fmincon} and non-linear semidefinite programming approaches converged to an infeasible point. Polynomial optimization exhibited convergence to the global optimum in three relaxations, Fig.~\ref{fig:girderC}. Notice, moreover, that the designs found for $r=\{1,2\}$ were of very high qualities.

\section{Conclusions}

In this contribution, we investigated four topology optimization techniques for the design of frame structures with a~given aspect-ratios all components of their cross-sections. Three of the optimization techniques---general non-linear formulation solved by the \texttt{fmincon} solver, optimality criteria, and non-linear semidefinite programming---provide a local solution to the considered optimization problem, and in turn converge to a~local optimum in general. Unfortunately, the NSDP and \texttt{fmincon} techniques also converged to infeasible points. From the local optimization approaches, the optimality criteria method seems to be the most efficient one. 

On the other hand, the last technique --- polynomial optimization --- allows to solve the problem to proven global optimality, while generating both lower and upper bounds in each of the relaxation order. Compared to the local approaches, the designs obtained in lower relaxation orders are of a~known quality with respect to the optimum. However, finding a proven global optimum requires a~considerably higher computational resources than solving the local optimization formulations. In all the test cases, the convergence of the PO hierarchy was finite.

\begin{nomenclature}
\item{\mathbf{0}}{Column vector or matrix of all zeros}
\item{\mathbf{1}}{Column vector or matrix of all ones}
\item{\mathbf{a}}{Cross-sectional areas column vector}
\item{\mathbf{a}_\mathrm{sc}}{Scaled cross-sectional areas column vector}
\item{a_i}{Cross-sectional area of the element $i$}
\item{a_{\mathrm{sc},i}}{Scaled cross-sectional area of the element $i$}
\item{\hat{\mathbf{a}}_r}{Optimized cross-sectional areas at the relaxation $r$}
\item{b_i}{Auxiliary variable associated with the element $i$}
\item{\mathbf{b}_r}{Polynomial space basis of maximum degree $r$}
\item{c}{Compliance (external work)}
\item{\tilde{c}}{Lower bound on compliance}
\item{\hat{c}}{Upper bound on compliance}
\item{c^*}{Globally optimal compliance}
\item{c_\mathrm{sc}}{Scaled compliance}
\item{c_{\mathrm{I},i}}{Positive constant}
\item{d_j}{Degree of the polynomial $j$}
\item{E}{Young modulus}
\item{\mathbf{f}}{Generalized force column vector}
\item{\mathbf{f}_i}{Generalized force column vector of element $i$}
\item{\hat{\mathbf{f}}}{Generalized force column vector associated with $\hat{\mathbf{K}}$}
\item{i}{Element index}
\item{j}{Polynomial index}
\item{I_i}{Moment of inertia of the element $i$}
\item{\mathbf{I}}{Identity matrix}
\item{k}{Iteration number}
\item{\mathbf{K}}{Stiffness matrix}
\item{\hat{\mathbf{K}}}{Positive definite principal submatrix of $\mathbf{K}$}
\item{\mathbf{K}_i}{Stiffness matrix of the element $i$}
\item{\mathcal{L}}{Lagrangian function}
\item{\bm{\ell}}{Column vector of element lengths}
\item{\ell_i}{Length of the element $i$}
\item{m}{Index of $\mathbf{x}$}
\item{\mathbf{M}_r}{Moment matrix of the order $r$}
\item{\mathbf{M}_{r-d_j}}{Localization matrix of the order $r-d_j$}
\item{n_\mathrm{e}}{Number of elements}
\item{n_\mathrm{n}}{Number of nodes}
\item{p_j}{Polynomial inequality}
\item{P_j}{Polynomial matrix inequality}
\item{\mathbf{p}_{\alpha,j}}{Coefficients of linear combinations of monomials}
\item{r}{Relaxation order, maximum polynomial degree}
\item{s}{Inverse of $c$}
\item{t_\mathrm{p}}{Parameter, web and flanges thickness}
\item{t}{Time}
\item{\mathbf{u}}{Generalized displacement column vector}
\item{\mathbf{v}}{Column vector of coefficients of linear combination}
\item{V}{Volume}
\item{\overline{V}}{Volume upper bound}
\item{\mathbf{x}}{Design variables column vector}
\item{x^\alpha}{Monomial indexed by $\bm{\alpha}$}
\item{y_\beta}{Moment associated with the basis $\mathbf{b}_r$}
\item{\mathbf{y}}{Column vector of moments}
\item{\mathbf{y}^*}{Optimal column vector of moments}
\item{\mathbf{y}_r}{Column vector of moments at $r$-th relaxation}

\item{\bm{\alpha}}{Vector of integers}
\item{\beta}{$\mathbf{b}_r$ index}
\item{\eta}{Tuning parameter, $0.3$ in this study}
\item{\zeta}{Move limit, $0.2$ in this study}
\item{\varepsilon}{Small positive number, $10^{-6}$ in this study}
\item{\bm{\lambda}}{Lagrange multiplier column vector}
\item{\mu}{Lagrange multiplier}
\item{\bm{\nu}}{Lagrange multiplier column vector}
\end{nomenclature}

\begin{acknowledgements}
The work of Marek Tyburec was supported by the Grant Agency of the CTU in Prague,  SGS19/033/OHK1/1T/11. Martin Kru\v{z}\'{i}k and Jan Zeman acknowledge the financial support by the European Regional Development Fund through the project CZ.02.1.01/0.0/0.0/16-019/0000778.
\end{acknowledgements}

\bibliographystyle{actapoly}
\bibliography{biblio}


\end{document}